\documentclass[12pt,reqno]{article}

\usepackage{enumerate}
\usepackage{enumitem}

\usepackage[usenames]{color}
\usepackage{amssymb}
\usepackage{amsmath}
\usepackage{amsthm}
\usepackage{amsfonts}
\usepackage{amscd}
\usepackage{graphicx}
\usepackage{enumitem}

\usepackage[colorlinks=true,
linkcolor=webgreen,
filecolor=webbrown,
citecolor=webgreen]{hyperref}

\definecolor{webgreen}{rgb}{0,.5,0}
\definecolor{webbrown}{rgb}{.6,0,0}

\usepackage{color}
\usepackage{fullpage}
\usepackage{float}

\usepackage{graphics}
\usepackage{latexsym}
\usepackage{epsf}
\usepackage{breakurl}

\theoremstyle{plain}
\newtheorem{thm}{Theorem}

\newtheorem{coro}[thm]{Corollary}

\theoremstyle{remark}

\theoremstyle{definition}

\numberwithin{equation}{section}

\begin{document}

\begin{center}
\vskip 1cm{\LARGE\bf
On the Number of  Words of Length $N$ over a Three-Letter Alphabet
\vskip .1in } \vskip 1cm

\large
Pablo Serra\\
Facultad de Matem\'atica, Astronom\'ia, F\'isica y Computaci\'on \\
Universidad Nacional de C\'ordoba, Universidad Nacional de C\'ordoba
 and
 Instituto de F\'isica Enrique Gaviola, CONICET\\
Av. Medina Allende S/N, CP:X5000HUA C\'ordoba,
 Argentina\\
\href{mailto:pablo.serra@unc.edu.ar}{\tt pablo.serra@unc.edu.ar}\\
\end{center}

\vskip .2in

\begin{abstract}
We study the classical problem of counting words formed from a three-letter alphabet by classifying them according to the residues modulo three of both the total length and the multiplicities of each letter. This approach yields explicit formulas for all cases and, as a consequence, several identities involving sums of trinomial coefficients.\\

\end{abstract}

\section{Problem statement}

The number of $N$-letter words that can be formed from an $m$-letter alphabet is a classical problem in combinatorics \textcolor{green}{\cite{stanley2012}}. For a configuration in which each letter $a_i$ appears $n_i \geq 0$, for $i=1,\ldots,m$, the number of distinct words is given by
\[
c(m;\{n_i\}) = \binom{N}{n_1,\ldots,n_m} ,  \text{ where } N = \sum_{i=1}^m n_i ,
\]
and
\[
\binom{N}{n_1,\ldots,n_m} =
\frac{(n_1+n_2+\cdots+n_m)!}{n_1!\,n_2!\cdots n_m!}
\]
denotes the multinomial coefficient. Summing over all configurations, the total number of $N$-letter words is
\begin{equation}
\label{ecgntc}
C(m;N) = \sum_{n_1+\cdots+n_m=N} c(m;\{n_i\})
= \sum_{n_1+\cdots+n_m=N} \binom{N}{n_1,\ldots,n_m}
= m^N.
\end{equation}

In this work, we focus on the case $m=3$, corresponding to an alphabet with three
letters $a_1, a_2, a_3$. We classify $N$-letter words according to the residue 
classes modulo $3$ of both $N$ and the multiplicities $n_i$, $i=1,2,3$. 
For convenience, we define
\[
n = N/3 ,
\]
and introduce $C_{\alpha\beta\gamma}(n)$, where $2 \geq \alpha \geq \beta \geq \gamma \geq 0$, as the number of $N$-letter words such that
\[
(n_1,n_2,n_3) \equiv (\alpha,\beta,\gamma) \pmod{3},
\]
including all permutations of $(n_1,n_2,n_3)$.

The possible values of $(\alpha,\beta,\gamma)$ depend on $N \bmod 3$:
\begin{itemize}
\item If $3 \mid N$: $(0,0,0)$, $(1,1,1)$, $(2,2,2)$, $(2,1,0)$;
\item If $N \equiv 1 \pmod{3}$: $(1,0,0)$, $(2,1,1)$, $(2,2,0)$;
\item If $N \equiv 2 \pmod{3}$: $(2,0,0)$, $(1,1,0)$, $(2,2,1)$.
\end{itemize}

Our main result is the following.

\begin{thm}
For $n \in \mathbb{N}$, the number of words of length $N=3n$ over a three-letter alphabet in the cases described above is given by
\begin{subequations}
\label{ect}
\begin{align}
\label{ecAt}
&C_{000}(n)=3^{3n-2} + (1+(-1)^n)\, i^n\, 3^{(3n-2)/2} ,\\
\label{ecBt}
&C_{111}(n)=3^{3n-2} - \left((1+(-1)^n)
+i\sqrt{3} (1-(-1)^n)\right)\frac{i^n}{2}\, 3^{(3n-2)/2} ,\\
\label{ecCt}
&C_{222}(n)=3^{3n-2} - \left((1+(-1)^n)
- i\sqrt{3} (1-(-1)^n)\right)\frac{i^n}{2} 3^{(3n-2)/2} \,,\\
\label{ecDt}
&C_{210}(n)=2 \cdot 3^{3n-1} .
\end{align}
\end{subequations}

The number of words of length $N=3n+1$ is
\begin{subequations}
\label{ec1}
\begin{align}
\label{ec100}
&C_{100}(n)=3^{3n}+ (1+(-1)^n) i^n 3^{3n/2} ,\\
\label{ec211}
&C_{211}(n)=3^{3n}- \left((1+(-1)^n)+ (1-(-1)^n)\sqrt{3} i\right)
\frac{i^n 3^{3n/2}}{2} ,\\
\label{ec220}
&C_{220}(n)=3^{3n}- \left((1+(-1)^n) -(1-(-1)^n)\sqrt{3} i\right)
\frac{i^n 3^{3n/2}}{2} .
\end{align}
\end{subequations}

The number of words of length $N=3n+2$ is
\begin{subequations}
\begin{align}
\label{ec200}
&C_{200}(n)=3^{3n+1}+ (1-(-1)^n)\, i^{n+1}\, 3^{(3n+1)/2} ,\\
\label{ec110}
&C_{110}(n)=3^{3n+1}+ \left((1+(-1)^n)\sqrt{3}- i(1-(-1)^n)\right)
\frac{i^n\, 3^{(3n+1)/2}}{2} ,\\
\label{ec221}
&C_{221}(n)=3^{3n+1}- \left((1+(-1)^n)\sqrt{3}+ i(1-(-1)^n)\right)
\frac{i^n\, 3^{(3n+1)/2}}{2} .
\end{align}
\end{subequations}
\end{thm}

The sequences in Eqs.~\eqref{ecAt}--\eqref{ecCt} have been included in the OEIS as A391468, A391469, and A391470, respectively \cite{A391468,A391469,A391470}. 
From Eq.~\eqref{ecgntc}, one must have
\begin{align*}
C_{000}(n)+C_{111}(n)+C_{222}(n)+C_{210}(n)&=3^{3n}, \\
C_{100}(n)+C_{211}(n)+C_{220}(n)&=3^{3n+1}, \\
C_{200}(n)+C_{110}(n)+C_{221}(n)&=3^{3n+2},
\end{align*}
which can be verified directly.

Notably, this theorem, together with Eq.~\eqref{ect}, provides new results for sums of trinomial coefficients:
\[
p_{\alpha\beta\gamma} \sum_{k_1+k_2+k_3=n}
\binom{N}{3k_1+\alpha,3k_2+\beta,3k_3+\gamma}
= C_{\alpha\beta\gamma}(n),
\]
where $p_{\alpha\beta\gamma}$ denotes the number of distinct permutations of the indices $\alpha$, $\beta$, and $\gamma$.

Since the corresponding generating functions are geometric series, the following 
corollary follows directly.
\begin{coro}
The generating functions of $C_{\alpha \beta \gamma}$, defined by

\[
g_{\alpha \beta \gamma}(x)=\sum_{n=0}^\infty C_{\alpha \beta \gamma}(n)x^n,
\]

\noindent are given by

\begin{align*}{ll}
&g_{000}(x)=\frac{1}{9(1-27\,x)}+\frac{2 (4+27 x^2)}{9(1+27 x^2)} ;
&g_{111}(x)=\frac{1}{9(1-27\,x)}+\frac{54 x^2+27 x-1}{9(1+27 \,x^2)};  \\
&g_{222}(x)=\frac{1}{9(1-27\,x)}+\frac{54 x^2-27 x-1}{9(1+27 \,x^2)} ; 
&g_{210}(x)=\frac{18\,x}{1-27\,x};\\
&g_{100}(x)=\frac{1}{1-27 x}+\frac{2}{1+27 x^2} ;
&g_{211}(x)=\frac{1}{1-27 x}+\frac{(9 x-1)}{1+27 x^2} ;\\
&g_{220}(x)=\frac{1}{1-27 x}- \frac{(9 x+1)}{1+27 x^2} ;
&g_{200}(x)=\frac{3}{1-27 x}- \frac{18 x}{1+27 x^2} ;\\
&g_{110}(x)=\frac{3}{1-27 x}+ \frac{3 (1+3 x)}{1+27 x^2} ;
&g_{221}(x)=\frac{3}{1-27 x}- \frac{3 (1-3 x)}{1+27 x^2} .
\end{align*}

\end{coro}

To prove the theorem, we first derive recurrence relations (RR) for each $C_{\alpha\beta\gamma}$ in the case $3 \mid N$, and then establish them by mathematical induction. Once these relations are proven, all cases with $N \not\equiv 0 \pmod{3}$ are completely determined by the corresponding solutions for $3 \mid N$.

\section{The case $3 \mid N$}
\subsection{Recurrence relations}

To determine the recurrence relations (RR) for the case $3 \mid N$, we make use of a well-known recurrence for multinomial coefficients \cite{ha1971}, which in the case of trinomial coefficients takes the form
\begin{equation}
\label{erct}
\binom{P}{p_1,\,p_2,\,p_3} =
\sum_{j=1}^3 \binom{P-1}{p_1-\delta_{1,j},\,p_2-\delta_{2,j},\,p_3-\delta_{3,j}},
 \text{with } P=p_1+p_2+p_3 ,
\end{equation}
where $\delta_{i,j}$ denotes the Kronecker delta,
\[
\delta_{i,j} =
\begin{cases}
1, & \text{if } i=j,\\
0, & \text{if } i\ne j.
\end{cases}
\]

In our setting, $3 \mid P$, and therefore Eq.~\eqref{erct} must be applied three times in order to derive the desired recurrence relations. Starting with $C_{000}(n)$, we obtain
\begin{align*}
&C_{000}(n)
= \sum_{k_1+k_2+k_3=n}
\binom{3(k_1+k_2+k_3)}{3k_1,\,3k_2,\,3k_3} \nonumber\\
&= \sum_{k_1+k_2+k_3=n}
\sum_{j_1,j_2,j_3=1}^3
\binom{3(k_1+k_2+k_3-1)}{
3k_1-\delta_{j_1,1}-\delta_{j_2,1}-\delta_{j_3,1},\,
3k_2-\delta_{j_1,2}-\delta_{j_2,2}-\delta_{j_3,2},\,
3k_3-\delta_{j_1,3}-\delta_{j_2,3}-\delta_{j_3,3}} \,.
\end{align*}

Using the permutation symmetry of the trinomial coefficients, this expression simplifies to
\begin{align*}
C_{000}(n)=\sum_{k_1+k_2+k_3=n}&\biggl(
3\,\binom{3(k_1+k_2+k_3-1)}{3(k_1-1),\,3k_2,\,3k_3} \\
&+18\,\binom{3(k_1+k_2+k_3-1)}{3k_1,\,3k_2-1,\,3k_3-2} \\
&+6\,\binom{3(k_1+k_2+k_3-1)}{3k_1-1,\,3k_2-1,\,3k_3-1}
\biggr),
\end{align*}
where we adopt the convention that any trinomial coefficient with a negative entry is zero. Using the definitions of $C_{\alpha\beta\gamma}$, we then obtain
\begin{equation}
\label{erpa}
C_{000}(n)=3\,C_{000}(n-1)+6\,C_{222}(n-1)+3\,C_{210}(n-1)\,.
\end{equation}

Proceeding analogously, we obtain the recurrence relations for the remaining cases. For $C_{111}(n)$,
\begin{equation}
\label{erpb}
C_{111}(n)=3\,C_{111}(n-1)+6\,C_{000}(n-1)+3\,C_{210}(n-1)\,,
\end{equation}
for $C_{222}(n)$,
\begin{equation}
\label{erpc}
C_{222}(n)=3\,C_{222}(n-1)+6\,C_{111}(n-1)+3\,C_{210}(n-1)\,,
\end{equation}
and for $C_{210}(n)$,
\begin{equation}
\label{erpd}
C_{210}(n)=18\Bigl(C_{210}(n-1)+C_{000}(n-1)+C_{111}(n-1)+C_{222}(n-1)\Bigr)\,.
\end{equation}

The next step is to decouple the recurrence relations given in Eqs.~\eqref{erpa}--\eqref{erpd}.

\subsection{Decoupling the recurrence relations}

Equation~\eqref{erpb} can be rewritten as
\begin{equation}
\label{erpbd}
C_{210}(n-1)=\frac{1}{3}\Bigl(C_{111}(n)-3C_{111}(n-1)-6C_{000}(n-1)\Bigr)\,.
\end{equation}

Substituting Eq.~\eqref{erpbd} into Eq.~\eqref{erpd} yields
\begin{equation}
\label{erpdc}
C_{222}(n)=\frac{1}{54}\Bigl(-6C_{000}(n+1)+54C_{000}(n)-21C_{111}(n+1)
+C_{111}(n+2)\Bigr)\,.
\end{equation}

Similarly, substituting Eq.~\eqref{erpbd} into Eq.~\eqref{erpc} gives
\begin{equation}
\label{erpcc}
6C_{000}(n-1)-C_{111}(n)-3C_{111}(n-1)+C_{222}(n)-3C_{222}(n-1)=0\,.
\end{equation}

From Eqs.~\eqref{erpdc} and \eqref{erpcc}, we obtain a relation involving only $C_{000}$ and $C_{111}$:
\begin{align}
\label{erpab1}
6C_{000}(n+1)-72C_{000}(n)-162C_{000}(n-1)
-C_{111}(n+2)+ \nonumber \\24C_{111}(n+1)-9C_{111}(n)+162C_{111}(n-1)=0\,.
\end{align}

On the other hand, substituting Eq.~\eqref{erpbd} into Eq.~\eqref{erpa} gives
\begin{align}
\label{erpcc2}
C_{000}(n)+3C_{000}(n-1)-C_{111}(n)+ \nonumber \\
3C_{111}(n-1)-6C_{222}(n-1)=0 ,
\end{align}
from which we eliminate $C_{222}$ to obtain
\begin{equation}
\label{erpab2}
15C_{000}(n)-27C_{000}(n-1)-C_{111}(n+1)+12C_{111}(n)+27C_{111}(n-1)=0\,.
\end{equation}

The coupled recurrence relations \eqref{erpab1} and \eqref{erpab2} can now be decoupled. Both $C_{000}$ and $C_{111}$ satisfy the same third-order recurrence relation, differing only in their initial conditions:
\begin{align}
\label{errsa}
C_{000}(n)=27\Bigl(C_{000}(n-1)-C_{000}(n-2)+27C_{000}(n-3)\Bigr), \nonumber\\
 C_{000}(1)=3,;C_{000}(2)=63,;C_{000}(3)=2187,
\end{align}
\begin{align}
\label{errsb}
C_{111}(n)=27\Bigl(C_{111}(n-1)-C_{111}(n-2)+27C_{111}(n-3)\Bigr), \nonumber\\
 C_{111}(1)=6,;C_{111}(2)=90,;C_{111}(3)=2106.
\end{align}

We now derive recurrence relations for $C_{222}$ and $C_{210}$. Subtracting three times Eq.~\eqref{erpc} from Eq.~\eqref{erpd}, we obtain
\[
C_{210}(n)-3C_{222}(n)=9\Bigl(2C_{000}(n-1)+C_{222}(n-1)+C_{210}(n-1)\Bigr).
\]

This identity allows us to eliminate $C_{000}$ from Eq.~\eqref{erpa}, yielding
\begin{equation}
\label{ecd1}
3C_{222}(n+1)+81C_{222}(n-1)-C_{210}(n+1)+12C_{210}(n)+27C_{210}(n-1)=0.
\end{equation}

A second relation between $C_{222}$ and $C_{210}$ is obtained by combining Eqs.~\eqref{erpd}, \eqref{erpb}, and \eqref{erpc}:
\begin{equation}
\label{esinb2}
C_{222}(n+1)+27C_{222}(n-1)-5C_{210}(n)+9C_{210}(n-1)=0.
\end{equation}

From Eqs.~\eqref{ecd1} and \eqref{esinb2}, we derive a recurrence relation for $C_{210}$:
\begin{equation}
\label{erpdf2}
C_{210}(n)=27 C_{210}(n-1),  C_{210}(1)=18,
\end{equation}
and, consequently, a recurrence relation for $C_{222}$:
\begin{align}
\label{err222}
C_{222}(n)=27\Bigl(C_{222}(n-1)-C_{222}(n-2)+27C_{222}(n-3)\Bigr), \nonumber\\
C_{222}(1)=0,;C_{222}(2)=90,;C_{222}(3)=2268.
\end{align}

This recurrence is identical to those satisfied by $C_{000}$ and $C_{111}$, differing only in the initial conditions. The fact that these three recurrence relations coincide can be understood in terms of symmetry. Consider the addition of three letters to a configuration of length $3n$:
\begin{itemize}
\item If all three letters are identical, each class maps to itself;
\item If all three letters are distinct, $(000)\to(111)\to(222)\to(000)$, while $(210)$ remains invariant;
\item If one letter is absent, $(000)$, $(111)$, and $(222)$ map to $(210)$, while $(210)$ can map to any class.
\end{itemize}

\subsection{Proof of the proposed solutions for the decoupled recurrence relations} 

In this section, we complete the proof of the theorem by mathematical induction for each expression of $C_{\alpha\beta\gamma}$.

\subsubsection{Proof of the proposed solution for $C_{000}$}

We show that Eq.~\eqref{ecAt} satisfies the recurrence relation for $C_{000}$ given in Eq.~\eqref{errsa}.

\begin{itemize}

\item \textbf{Verification of a specific case.} We evaluate $C_{000}(4)$ using both Eq.~\eqref{errsa} and Eq.~\eqref{ecAt}:
\begin{align*}
C_{000}(4)
&=27\left(C_{000}(3)-C_{000}(2)+27\,C_{000}(1)\right) \nonumber\\
&=27(2187-63+27\cdot 3)=59535 \,, \nonumber\\
C_{000}(4)
&=3^{3\cdot4-2}+\left(1+(-1)^4\right)i^4\,3^{(3\cdot4-2)/2} \nonumber\\
&=3^{10}+2\cdot 3^5=59535 \,.
\end{align*}

\item \textbf{Inductive hypothesis.} Assume that Eq.~\eqref{ecAt} holds for all $n=1,\ldots,k$.

\item \textbf{Inductive step.} For $n=k+1$, using the recurrence relation and the inductive hypothesis (IH), we obtain
\begin{align*}
C_{000}(k+1)
&=27\Bigl(C_{000}(k)-C_{000}(k-1)+27\,C_{000}(k-2)\Bigr) \nonumber\\
&\stackrel{\mathrm{IH}}{=}27\Bigl(
3^{3k-2}+(1+(-1)^k)i^k\,3^{(3k-2)/2} \nonumber\\
& -\bigl(3^{3(k-1)-2}+(1+(-1)^{k-1})i^{k-1}\,3^{(3(k-1)-2)/2}\bigr) \nonumber\\
& +27\bigl(3^{3(k-2)-2}+(1+(-1)^{k-2})i^{k-2}\,3^{(3(k-2)-2)/2}\bigr)
\Bigr) \nonumber\\
&=3^{3(k+1)-2}+(1+(-1)^{k+1})i^{k+1}\,3^{(3(k+1)-2)/2} \,.
\end{align*}

\end{itemize}

This completes the proof by induction. 

\subsubsection{Proof of the proposed solution for $C_{111}$}

We show that Eq.~\eqref{ecBt} satisfies the recurrence relation for $C_{111}$ given in Eq.~\eqref{errsb}, thereby proving that it is the solution (item $B$ of the theorem).

\begin{itemize}

\item \textbf{Verification of a specific case.} Evaluating both Eq.~\eqref{errsb} and Eq.~\eqref{ecBt} at $n=4$, we obtain
\begin{align*}
C_{111}(4)
&=27\bigl(C_{111}(3)-C_{111}(2)+27\,C_{111}(1)\bigr) \nonumber\\
&=27(2106-90+27\cdot 6)=58806 \,, \nonumber\\
C_{111}(4)
&=3^{3\cdot4-2}
-\Bigl((1+(-1)^4)+(1-(-1)^4)i\sqrt{3}\Bigr)\frac{i^4}{2}\,3^{(3\cdot4-2)/2} \nonumber\\
&=3^{10}-3^5=58806 \,.
\end{align*}

\item \textbf{Inductive hypothesis.} Assume that Eq.~\eqref{ecBt} holds for all $n=1,\ldots,k$.

\item \textbf{Inductive step.} For $n=k+1$, using the recurrence relation and the inductive hypothesis (IH), we obtain
\begin{align*}
C_{111}(k+1)
&=27\Bigl(C_{111}(k)-C_{111}(k-1)+27\,C_{111}(k-2)\Bigr) \nonumber\\
&\stackrel{\mathrm{IH}}{=}27\Bigl(
3^{3k-2}
-\Bigl((1+(-1)^k)+(1-(-1)^k)i\sqrt{3}\Bigr)\frac{i^k}{2}\,3^{(3k-2)/2}
\nonumber\\
& -\Bigl(3^{3(k-1)-2}
-\Bigl((1+(-1)^{k-1})+(1-(-1)^{k-1})i\sqrt{3}\Bigr)
\frac{i^{k-1}}{2}\,3^{(3(k-1)-2)/2}\Bigr)
\nonumber\\
& +27\Bigl(3^{3(k-2)-2}
-\Bigl((1+(-1)^{k-2})+(1-(-1)^{k-2})i\sqrt{3}\Bigr)
\frac{i^{k-2}}{2}\,3^{(3(k-2)-2)/2}\Bigr)
\Bigr) \nonumber\\
&=3^{3(k+1)-2}
-\Bigl((1+(-1)^{k+1})+(1-(-1)^{k+1})i\sqrt{3}\Bigr)
\frac{i^{k+1}}{2}\,3^{(3(k+1)-2)/2} \,.
\end{align*}

\end{itemize}

This completes the proof by induction. 

\subsubsection{Proof of the proposed solution for $C_{222}$}

Similarly to the previous cases, we prove that Eq.~\eqref{ecCt} is a solution of the recurrence relation for $C_{222}$ given in Eq.~\eqref{err222}, by mathematical induction.

\begin{itemize}

\item \textbf{Verification of a specific case.} Evaluating both Eq.~\eqref{err222}
 and Eq.~\eqref{ecCt} at $n=4$, we obtain
\begin{align*}
C_{222}(4)
&=27\bigl(C_{222}(3)-C_{222}(2)+27\,C_{222}(1)\bigr) \nonumber\\
&=27(2268-90)=58806 \,, \nonumber\\
C_{222}(4)
&=3^{3\cdot4-2}
-\Bigl((1+(-1)^4)-i\sqrt{3}(1-(-1)^4)\Bigr)\frac{i^4}{2}\,3^{(3\cdot4-2)/2} \nonumber\\
&=3^{10}-3^5=58806 \,.
\end{align*}

\item \textbf{Inductive hypothesis.} Assume that Eq.~\eqref{ecCt} holds for all $n=1,\ldots,k$.

\item \textbf{Inductive step.} For $n=k+1$, using the recurrence relation and the inductive hypothesis (IH), we obtain
\begin{align*}
C_{222}(k+1)
&=27\Bigl(C_{222}(k)-C_{222}(k-1)+27\,C_{222}(k-2)\Bigr) \nonumber\\
&\stackrel{\mathrm{IH}}{=}27\Bigl(
3^{3k-2}
-\Bigl((1+(-1)^k)-i\sqrt{3}(1-(-1)^k)\Bigr)\frac{i^k}{2}\,3^{(3k-2)/2}
\nonumber\\
& -\Bigl(3^{3(k-1)-2}
-\Bigl((1+(-1)^{k-1})-i\sqrt{3}(1-(-1)^{k-1})\Bigr)
\frac{i^{k-1}}{2}\,3^{(3(k-1)-2)/2}\Bigr)
\nonumber\\
& +27\Bigl(3^{3(k-2)-2}
-\Bigl((1+(-1)^{k-2})-i\sqrt{3}(1-(-1)^{k-2})\Bigr)
\frac{i^{k-2}}{2}\,3^{(3(k-2)-2)/2}\Bigr)
\Bigr) \nonumber\\
&=3^{3(k+1)-2}
-\Bigl((1+(-1)^{k+1})-i\sqrt{3}(1-(-1)^{k+1})\Bigr)
\frac{i^{k+1}}{2}\,3^{(3(k+1)-2)/2} \,.
\end{align*}

\end{itemize}

This completes the proof by induction. 

\subsubsection{Proof of the proposed solution for $C_{210}$}

Finally, to complete the demonstration of the theorem, we have to prove that
Eq.~\eqref{ecDt} is the solution of   Eq.~\eqref{erpdf2},

\begin{itemize}
\item $k=2:$
\[
        C_{210}(2) = 2\cdot 3^{3\cdot 2-1} = 2\cdot 3^5 = 486 ;
C_{210}(2) = 27 C_{210}(1) = 27 \cdot 18 =486 .
\]

\item Assume  valid the Eq.~\eqref{ecDt} for   $n=k$.

\item $n=k+1:$
\[
C_{210}(k+1)\,=\,3^3\,C_{210}(k)\,\stackrel{IH}{=}\,3^3\,2\cdot 3^{3\,k-1}\,=\,
2\cdot 3^{3k+2}\,=\, 2\cdot 3^{3(k+1)-1}. 
\]
\end{itemize}

which completes the proof of the first case of theorem.

The fact that $C_{000}, C_{111}$, and $C_{222}$ have the same recurrence 
relations, with different initial condition can be understood in terms of 
symmetry. Adding three letters to a configuration of length $3n$:
\begin{itemize}
\item If all three letters are identical, each class maps to itself;
\item If all three letters are distinct, $(000)\to(111)\to(222)\to(000)$, while $(210)$ remains invariant;
\item If one letter is absent, $(000)$, $(111)$, and $(222)$ map to $(210)$, while $(210)$ can map to any class.
\end{itemize}

Alternatively, the expressions given in Eqs.~\eqref{ecAt}--\eqref{ecDt} can be 
verified simultaneously by a unified argument based on recurrence relations and 
their characteristic polynomials.

The sequences $C_{000}(n)$, $C_{111}(n)$, and $C_{222}(n)$ satisfy the same third-order linear recurrence relation,
\begin{equation}
\label{errsunified}
C(n)=27\left(C(n-1)-C(n-2)+27\,C(n-3)\right),
\end{equation}
with different initial conditions. The characteristic polynomial associated with this recurrence is
\[
x^3 - 27(x^2 - x + 27),
\]
whose roots are
\[
x_1=27,  x_2=3^{3/2} i,  x_3=-3^{3/2} i.
\]

Therefore, the general solution of Eq.~\eqref{errsunified} is of the form
\[
C(n)=A\,x_1^n + B\,x_2^n + D\,x_3^n,
\]
where the constants $A,B,D$ are determined by the initial values of each sequence.

By matching the initial conditions for $C_{000}$, $C_{111}$, and $C_{222}$, we recover precisely the expressions given in Eqs.~\eqref{ecAt}, \eqref{ecBt}, and \eqref{ecCt}, respectively.

Similarly, the sequence $C_{210}(n)$ satisfies the first-order recurrence relation
\begin{equation}
C_{210}(n)=27\,C_{210}(n-1),
\end{equation}
whose general solution is
\[
C_{210}(n)=K\,27^n.
\]
Using the initial condition $C_{210}(1)=18$, we obtain
\[
C_{210}(n)=2\cdot 3^{3n-1},
\]
which coincides with Eq.~\eqref{ecDt}.

Finally, direct substitution shows that the expressions in 
Eqs.~\eqref{ecAt}--\eqref{ecDt} satisfy their respective recurrence relations 
and initial conditions. By uniqueness of solutions to linear recurrence 
relations, these expressions are the desired closed forms, giving an alternative 
proof of the first case, $3 \mid N$, of the theorem.

\section{Case $N \equiv 1 \mod(3)$}

For $N \equiv 1 \pmod{3}$, we define $n = (N-1)/3$. In this case, there are three possible configurations:
\begin{enumerate}
\item \label{ss100} $(n_1, n_2, n_3) \equiv (1,0,0) \pmod{3}$,
\item \label{ss211} $(n_1, n_2, n_3) \equiv (2,1,1) \pmod{3}$,
\item \label{ss220} $(n_1, n_2, n_3) \equiv (2,2,0) \pmod{3}$.
\end{enumerate}

\subsection{Case \ref{ss100}: $(n_1, n_2, n_3) \equiv (1,0,0) \pmod{3}$}

A word satisfying this constraint can be obtained from a word of length $N-1$ in two different ways:

\begin{enumerate}[label=\roman*.]
\item By appending an arbitrary letter to a word with $(n_1, n_2, n_3) \equiv (0,0,0) \pmod{3}$.
\item By appending a letter $a_i$ such that $n_i \equiv 2 \pmod{3}$ to a word with $(n_1, n_2, n_3) \equiv (2,1,0) \pmod{3}$.
\end{enumerate}

\noindent Consequently, $C_{100}(n)$ satisfies
\[
C_{100}(n)=3\,C_{000}(n)+C_{210}(n)
=3^{3n}+(1+(-1)^n)i^n\,3^{3n/2}\,.
\]

\subsection{Case \ref{ss211}: $(n_1, n_2, n_3) \equiv (2,1,1) \pmod{3}$}

A configuration of this type can be obtained from a word of length $N-1$ in two different ways:

\begin{enumerate}[label=\roman*.]
\item By appending an arbitrary letter to a word with $(n_1,n_2,n_3)\equiv (1,1,1)\pmod{3}$.
\item By appending a letter $a_i$ such that $3 \mid n_i$ to a word with $(n_1,n_2,n_3)\equiv (2,1,0)\pmod{3}$.
\end{enumerate}

\noindent Consequently,
\[
C_{211}(n)=3\,C_{111}(n)+C_{210}(n)
=3^{3n}-\Bigl((1+(-1)^n)+(1-(-1)^n)\sqrt{3}\,i\Bigr)\frac{i^n\,3^{3n/2}}{2}.
\]

\subsection{Case \ref{ss220}: $(n_1, n_2, n_3) \equiv (2,2,0) \pmod{3}$}

This configuration can be obtained from the case $3 \mid (N-1)$ in two ways:

\begin{enumerate}[label=\roman*.]
\item From $(n_1,n_2,n_3)\equiv (2,2,2)\pmod{3}$ by appending an arbitrary letter.
\item From $(n_1,n_2,n_3)\equiv (2,1,0)\pmod{3}$ by appending a letter $a_i$ such that $n_i \equiv 1 \pmod{3}$.
\end{enumerate}

\noindent Therefore,
\[
C_{220}(n)=3\,C_{222}(n)+C_{210}(n)
=3^{3n}-\Bigl((1+(-1)^n)-(1-(-1)^n)\sqrt{3}\,i\Bigr)\frac{i^n\,3^{3n/2}}{2}.
\]

\section{Case $N \equiv 2 \pmod{3}$}

If $N \equiv 2 \pmod{3}$, we define $n=(N-2)/3$. In this case, there are three possible configurations:
\begin{enumerate}
\item \label{ss200} $(n_1, n_2, n_3) \equiv (2,0,0) \pmod{3}$,
\item \label{ss110} $(n_1, n_2, n_3) \equiv (1,1,0) \pmod{3}$,
\item \label{ss221} $(n_1, n_2, n_3) \equiv (2,2,1) \pmod{3}$.
\end{enumerate}

\subsection{Case \ref{ss200}: $(n_1, n_2, n_3) \equiv (2,0,0) \pmod{3}$}

This configuration can be obtained from the case $N \equiv 1 \pmod{3}$ in two ways:

\begin{enumerate}[label=\roman*.]
\item From $(n_1,n_2,n_3)\equiv (1,0,0)\pmod{3}$ by appending a letter $a_1$.
\item From $(n_1,n_2,n_3)\equiv (2,2,0)\pmod{3}$ by appending a letter corresponding to an index $i$ such that $n_i \equiv 2 \pmod{3}$.
\end{enumerate}

\noindent Hence,
\[
C_{200}(n)=C_{100}(n)+2\,C_{220}(n)
=3^{3n+1}+(1-(-1)^n)i^{n+1}\,3^{(3n+1)/2}.
\]

\subsection{Case \ref{ss110}: $(n_1, n_2, n_3) \equiv (1,1,0) \pmod{3}$}

This configuration can be obtained from the case $N \equiv 1 \pmod{3}$ in two ways:

\begin{enumerate}[label=\roman*.]
\item From $(n_1,n_2,n_3)\equiv (1,0,0)\pmod{3}$ by appending either $a_2$ or $a_3$.
\item From $(n_1,n_2,n_3)\equiv (2,1,1)\pmod{3}$ by appending $a_1$.
\end{enumerate}

\noindent Therefore,
\[
C_{110}(n)=2\,C_{100}(n)+C_{211}(n)
=3^{3n+1}+\Bigl((1+(-1)^n)\sqrt{3}-i(1-(-1)^n)\Bigr)\frac{i^n\,3^{(3n+1)/2}}{2}.
\]

\subsection{Case \ref{ss221}: $(n_1, n_2, n_3) \equiv (2,2,1) \pmod{3}$}

This configuration can be obtained from the case $N \equiv 1 \pmod{3}$ in two ways:

\begin{enumerate}[label=\roman*.]
\item From $(n_1,n_2,n_3)\equiv (2,1,1)\pmod{3}$ by appending either $a_2$ or $a_3$.
\item From $(n_1,n_2,n_3)\equiv (2,2,0)\pmod{3}$ by appending a letter corresponding to the index with residue $0$.
\end{enumerate}

\noindent Hence,
\[
C_{221}(n)=2\,C_{211}(n)+C_{220}(n)
=3^{3n+1}-\Bigl((1+(-1)^n)\sqrt{3}+i(1-(-1)^n)\Bigr)\frac{i^n\,3^{(3n+1)/2}}{2}.
\]

This complete the proof of the theorem.

\section{Acknowledgments}

The author thanks Juan Pablo Rossetti for fruitful discussions and for 
pointing out the symmetry properties of the recurrence relations. The author 
also thanks the editors of the OEIS for helpful discussions. This work was 
supported in part by CONICET and SECYT-UNC.

\end{document}